\documentclass[graybox]{svmult}
\usepackage{mathptmx}
\usepackage{helvet}
\usepackage{courier}
\usepackage{type1cm}
\usepackage{makeidx}
\usepackage{graphicx}
\usepackage{multicol}
\usepackage[bottom]{footmisc}
\usepackage{amsmath, amssymb}
\usepackage[compatibility=false,size=footnotesize]{caption}
\usepackage{subcaption}
\usepackage{hyperref}
\usepackage{tikz}
\usetikzlibrary{shapes,patterns,snakes} 
\renewcommand{\div}{\nabla\cdot}
\newcommand{\curl}{\nabla\times}
\begin{document}
\title*{Exploring Parallel-in-Time Approaches for Eddy Current Problems}
\author{Stephanie Friedhoff, Jens Hahne, Iryna Kulchytska-Ruchka, and Sebastian Sch\"ops}
\institute{Stephanie Friedhoff and Jens Hahne \at Fakult\"at f\"ur Mathematik und Naturwissenschaften, Bergische Universit\"at Wuppertal, Gau{\ss}str. 20, 42119 Wuppertal, Germany, \email{friedhoff@math.uni-wuppertal.de} and \email{jens.hahne@math.uni-wuppertal.de}
\and Sebastian Sch\"ops and Iryna Kulchytska-Ruchka \at Centre for Computational Engineering, Technische Universit\"at Darmstadt, Dolivostr. 15, 64293 Darmstadt, Germany, \email{schoeps@temf.tu-darmstadt.de} and \email{kulchytska@temf.tu-darmstadt.de}}
\maketitle
\abstract*{}

\abstract{We consider the usage of parallel-in-time algorithms of the Parareal and multigrid-reduction-in-time (MGRIT) methodologies for the parallel-in-time solution of the eddy current problem. Via application of these methods to a two-dimensional model problem for a coaxial cable model, we show that a significant speedup can be achieved in comparison to sequential time stepping.}

\section{Introduction}
\label{sec:intro}
Recently, efficient and robust designs of electromechanical energy converters are gaining again in importance because of the transition towards sustainable energy in Europe (`Energiewende' in German). Electrical machinery is well understood and developed in industry close to their technical limits, but often without transient analysis or consideration of uncertainties in the design process. Such studies are only carried out late in the development process due to their high computational costs. This may lead to the fact that better or more robust designs are ruled out and not considered further on. One promising way to speed up transient analysis are parallel-in-time methods. 

In contrast to classical time-integration techniques based on a time-stepping approach, i.\,e., solving sequentially for one time step after the other, parallel-in-time algorithms allow simultaneous solution across multiple time steps. Starting with the work of Nievergelt \cite{Nivergelt1964}, various approaches for parallel-in-time integration have been explored; a recent review of the extensive literature in this area is \cite{Gander2015_Review}. The key practical aspect for choosing one of the many time-parallel methods when aiming at adding parallelism to an existing application code is the level of intrusiveness, i.\,e., the required amount of implementation effort. There are only a few time-parallel methods that are non-intrusive. In this paper, we consider two of these approaches, the Parareal method \cite{Lions_etal2001} and the multigrid-reduction-in-time (MGRIT) algorithm \cite{Falgout_etal_2014} that, in a specific two-level setting, can be viewed as a Parareal-type algorithm.
\section{Eddy current model problem}
\label{sec:model_problem}
\begin{figure}[t]
	\centering
	\begin{subfigure}[b]{0.7\textwidth}
		\centering
		\setlength{\unitlength}{2.9cm}
		\begin{picture}(7,0.5)%
			\put(0,0.1){\includegraphics[width=0.9\textwidth]{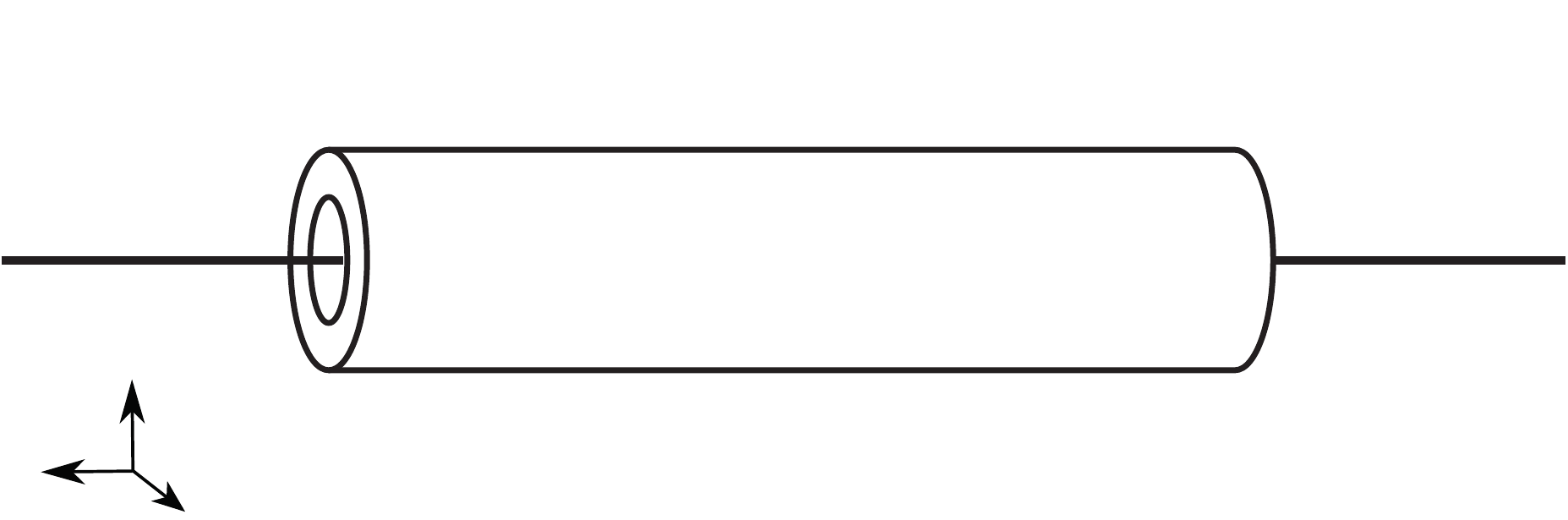}}%
			\put(0.31,0.09){\makebox(0,0)[lb]{\smash{$x$}}}%
			\put(0,0.15){\makebox(0,0)[lb]{\smash{$z$}}}%
			\put(0.2,0.37){\makebox(0,0)[lb]{\smash{$y$}}}%
		\end{picture}%
		\vspace{-2.5em}
		\caption{Cable model}
		\label{fig:wire}
	\end{subfigure}
	\begin{subfigure}[b]{0.25\textwidth}
		\centering
		\setlength{\unitlength}{2.9cm}
		\begin{picture}(1,1)%
			\put(-0.05,0){\includegraphics[width=\unitlength]{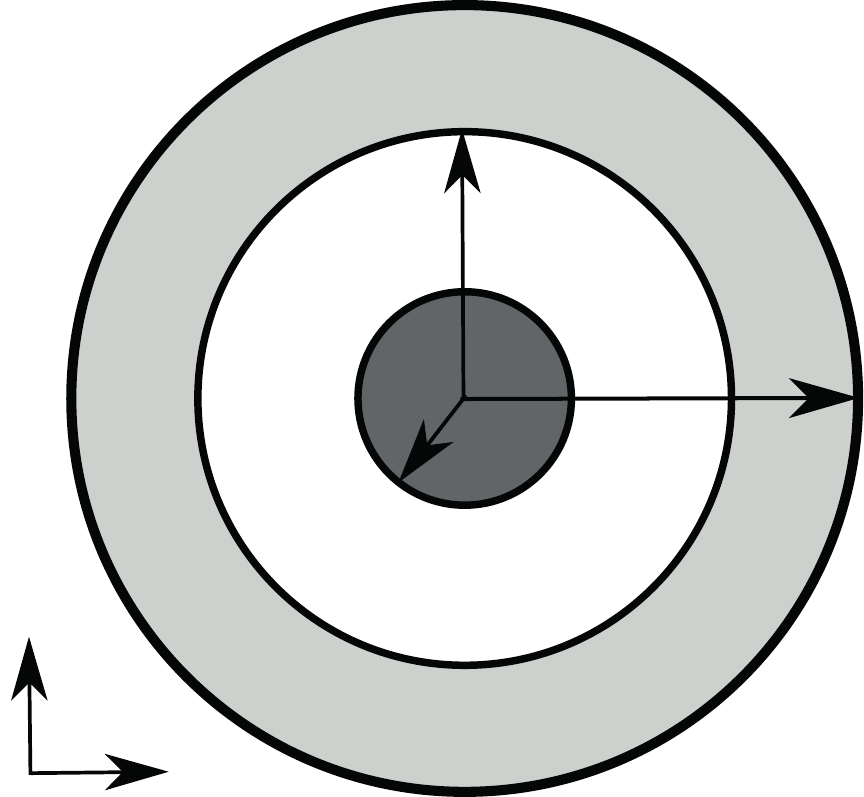}}%
			\put(0.82,0.52){\makebox(0,0)[lb]{\smash{$r_2$}}}%
			\put(0.5242171,0.65){\makebox(0,0)[lb]{\smash{$r_1$}}}%
			\put(0.40414256,0.48957478){\makebox(0,0)[lb]{\smash{$r_0$}}}%
			\put(0.17,0.01){\makebox(0,0)[lb]{\smash{$x$}}}%
			\put(-0.03,0.22){\makebox(0,0)[lb]{\smash{$y$}}}%
		\end{picture}%
		\caption{Cross-section}
		\label{fig:model}
	\end{subfigure}
	\caption{Cable model (a) and its cross section (b), dark grey region $\Omega_0$ models the wire, the white region an insulator $\Omega_1$ and the light grey annulus the conducting shield $\Omega_2$, see \cite{FEMM}}.\vspace{-1em}
\end{figure}
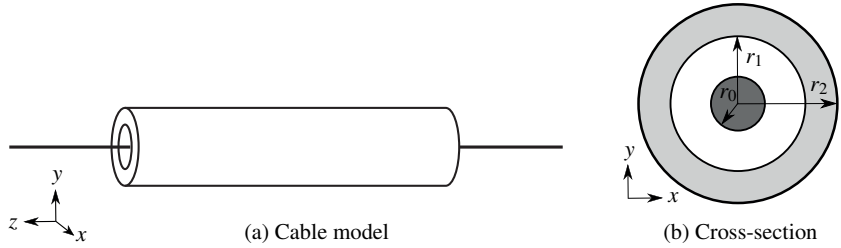

For an open, bounded domain $\Omega \subset \mathbb{R}^3$ and 
$t \in \mathcal{I} = (t_0,t_\mathrm{end}] \subset \mathbb{R}_{\geq0}$, the evolution 
of electromagnetic fields is governed by Maxwell's equations on $\Omega \times \mathcal{I}$ \cite{Jackson_1998aa}
\begin{align}
	\curl{\vec{E}}  &= -\partial_t \vec{B} ,
	&\quad
	\curl{\vec{H}}  &= \partial_t \vec{D}+\vec{J},
	&\quad
	\div{\vec{B}} &= 0,
	&\quad
	\div{\vec{D}} &= \rho,
	\label{eq:maxwell-equations}
\end{align}
with suitable initial and boundary conditions at time $t_0$ and $\partial\Omega$, respectively.
These equations are completed by constitutive relations
laws
\begin{align}
    \vec{D}=\varepsilon\vec{E}, \qquad 
    \vec{J}=\sigma \vec{E} + \vec{J}_\text{s}, \qquad 
    \vec{B}=\mu\vec{H}.
  \label{eq:const-laws-1}
\end{align}
In these equations, $\vec{H}$ is the magnetic field [A/m], $\vec{B}$ the magnetic 
flux density [T], $\vec{E}$ the electric field [V/m], $\vec{D}$ the electric flux 
density [C/m$^2$], $\vec{J}$ and $\vec{J}_\mathrm{s}$ 
are the total and source current density [A/m$^2$], $\rho$ 
is the electric charge density [C/m$^3$]. All fields are functions of space $\vec{x}\in\Omega$ and time $t\in\mathcal{I}$.
The material properties $\sigma\geq0$, $\varepsilon>0$ 
and $\mu>0$ are the electric conductivity, the electric permittivity and the magnetic 
permeability, respectively. It is convenient to invert the magnetic material law, i.\,e., 
$\vec{H}=\nu\vec{B}$, using the reluctivity $\nu$, where $\nu(B)$  can be a sufficiently smooth and bounded function of the magnitude $B=\|\vec{B}\|$, see \cite{Heise_1994aa}. In the following, we consider only devices where the displacement current can be neglected with respect to the source currents, i.\,e.,
$\|\partial_t \vec{D}\|\ll\|\vec{J}_\text{s}\|$.
An analysis of this error can be found in \cite{Schmidt_2008aa}. Assuming $\partial_t \vec{D}=0$ yields the so-called \emph{magnetoquasistatic} approximation or \emph{eddy current problem}. Eddy currents lead to the \emph{skin effect}, i.\,e., currents through a conductor are pushed to the surface if frequency increases \cite[Chapter 5.18]{Jackson_1998aa}.

One may introduce the (`modified' \cite{Emson_1988aa}) magnetic vector potential $\vec{A}$ such that $\vec{E}=-\partial_t\vec{A}$. Then, inserting the equations into each other yields
\begin{equation}
\sigma \partial_t \vec{A} + \nabla \times \big(\nu(\|\nabla\times\vec{A}\|) \nabla\times\vec{A}\big) = \vec{J}_\mathrm{s}.
\label{eq:pde}
\end{equation}
The source is defined as $\vec{J}_\mathrm{s}|_{\Omega_0}=\vec{e}_z/(\pi r_0^2) f_{n}(t)$ in the inner cable $\Omega_0$ and vanishes elsewhere, $\vec{e}_z$ denotes the unit vector in $z$-direction and the excitation is given by
\begin{equation}\label{eq:input_pwm}
	f_{n}(t)=\begin{cases}
	\mathrm{sign}\left[\sin\left(\dfrac{2\pi}{T} t\right)\right],\ & s_n(t)-\left|\sin\left(\dfrac{2\pi}{T}t\right)\right|<0,\\
	0,\ & \mathrm{otherwise,}
	\end{cases}
\end{equation}
where $s_n(t)=n/Tt -\left\lfloor n/Tt\right\rfloor,$ $t\in[0,T]$ is the common sawtooth pattern, with \linebreak$n=1100$ teeth and period $T=0.02\;$s \cite{Kulchytska-Ruchka_2018ac}. The reluctivity $\nu$ is modeled as vacuum ($1/\mu_0$) in $\Omega_0$ and $\Omega_1$, and is given in $\Omega_2$ by a Spline curve, the conductivity $\sigma$ is only non-zero in the tube region {$\Omega_2$} ($10$ MS/m).

Finally, the Ritz-Galerkin approach is employed with lowest order ansatz functions in space. When considering planar 2D problems, edge shape functions only have a $z$-component and can be constructed from the nodal shape functions $N_i(\vec{x})$ as
\begin{equation}
	\vec{A}=\sum_{i=1}^{N_\mathrm{dof}}\mathbf{u}_i \vec{w}_i(\vec{x})
	\quad\text{ with }\quad
	\vec{w}_i(\vec{x})=\frac{N_i(\vec{x})}{l_z}\vec{e}_z,
\end{equation}
where $l_z$ refers to the length in $z$-direction, $N_\mathrm{dof}=2269$. This leads to the equation 
\begin{equation}
\mathbf{M}_\sigma \mathbf{u}'+\mathbf{K}_\nu(\mathbf{u})\mathbf{u}=\mathbf{j}_\mathrm{s},
\label{eq:magnetoquasi_discret}
\end{equation}
with the matrices and the right-hand side
\begin{equation*}
M_{\sigma,i,j}=\!\!\int_{\Omega} \sigma\vec{w}_j\cdot\vec{w}_i\;\mathrm{d}\vec{x},
\text{~}
K_{\nu,i,j}(\cdot)=\!\!\int_{\Omega}\!\!\nu(\cdot)\nabla\times \vec{w}_j\cdot\nabla\times \vec{w}_i\;\mathrm{d}\vec{x}, 
\text{~}
j_{\mathrm{s},i}= \int_{\Omega}\!\!\vec{J}_\mathrm{s}\cdot\vec{w}_i\;\mathrm{d}\vec{x},
\end{equation*}
respectively. The resulting system \eqref{eq:magnetoquasi_discret} consists of differential-algebraic equations of index-1 due the vanishing entries $M_{\sigma,i,j}$ of the mass matrix $\mathbf{M}_\sigma$ in $\Omega_0$ and $\Omega_1$, \cite{Nicolet_1996aa}.

\section{Multigrid reduction in time}
\label{sec:mgrit}
The multigrid-reduction-in-time (MGRIT) algorithm \cite{Falgout_etal_2014} is an iterative, parallel method,
based on applying multigrid reduction (MGR) \cite{MRies_UTrottenberg_1979,MRies_etal_1983} principles in time, 
for solving time-stepping problems of the form
\begin{equation}\label{eq:ode:system}
	\mathbf{u}'(t) = \mathbf{f}(t,\mathbf{u}(t)), \quad \mathbf{u}(t_0) = \mathbf{g}_0, \quad t\in(t_0,t_\mathrm{end}] \subset \mathbb{R}_{\geq0},
\end{equation}
with initial condition, $\mathbf{g}_0$, at $t=t_0$. Note that form \eqref{eq:ode:system} can be a system of ODEs, arising, for example, after spatial discretization of a space-time PDE, or it can be a system of DAEs such as given in Equation \eqref{eq:magnetoquasi_discret}. Discretizing the time interval on a grid $t_i = i\Delta t$, $i=0,1,\ldots,N_t$, with, for notational convenience, constant time step $\Delta t=\left(t_\mathrm{end}-t_0\right)/N_t>0$, let $\mathbf{u}_i$ be an approximation to $\mathbf{u}(t_i)$ for $i=1,\ldots,N_t$, and let $\mathbf{u}_0 = \mathbf{u}(t_0)$. Then, considering a one-step time-independent time integration method with time-stepping operator, $\Phi_{\Delta t}$, that takes a solution at time $t_{i-1}$ to that at time $t_i$, along with a time-dependent forcing term, $\mathbf{g}_i$, the solution to~\eqref{eq:ode:system} is defined via time-stepping, which can also be represented as a forward solve of the linear system, written in block form as
\begin{equation}\label{eq:fine_linear_system}
	\mathbf{A}\mathbf{u}\equiv\begin{bmatrix}
		I\\
		-\Phi_{\Delta t} & I\\
		& \ddots & \ddots\\
		& & -\Phi_{\Delta t} & I
	\end{bmatrix}\begin{bmatrix}
		\mathbf{u}_0\\
		\mathbf{u}_1\\
		\vdots\\
		\mathbf{u}_{N_t}
	\end{bmatrix} = \begin{bmatrix}
		\mathbf{g}_0\\
		\mathbf{g}_1\\
		\vdots\\
		\mathbf{g}_{N_t}
	\end{bmatrix}\equiv \mathbf{g}.
\end{equation}
Hence, in the time dimension, this forward solve is completely sequential. 

Alternatively, considering the lower block bidiagonal structure, we could apply cyclic reduction, whereby we first solve the Schur complement system,
\begin{equation}\label{eq:sc_linear_system}
	\mathbf{A}_S\mathbf{u}_\Delta \equiv \begin{bmatrix}
		I\\
		-\Phi_{\Delta t}^m & I\\
		& \ddots & \ddots\\
		& & -\Phi_{\Delta t}^m & I
	\end{bmatrix} \begin{bmatrix}
		\mathbf{u}_0\\
		\mathbf{u}_m\\
		\vdots\\
		\mathbf{u}_{N_t}
	\end{bmatrix} = \begin{bmatrix}
		\mathbf{g}_0\\
		\hat{\mathbf{g}}_m\\
		\vdots\\
		\hat{\mathbf{g}}_{N_t}
	\end{bmatrix} \equiv \hat{\mathbf{g}},
\end{equation}
for the value of the solution at every $m$-th temporal point, with consistently restricted forcing terms. Then define the solution at the remaining temporal points by local (and parallel) time-stepping between those points defined from the Schur complement. MGRIT is based on interpreting this cyclic reduction approach as a two-level MGR algorithm, enabling parallelism in the solution process \eqref{eq:fine_linear_system}. Therefore, define a coarse temporal mesh, or (using multigrid terminology) the set of C-points, to be those points included in the Schur complement system \eqref{eq:sc_linear_system}, with the remaining temporal points as the set of F-points. Further define ``ideal'' interpolation as the map which takes the solution at the C-points and yields a zero residual at the F-points, with a similar definition for ``ideal'' restriction. The Schur complement then arises as the standard Petrov-Galerkin coarse-grid operator with these definitions of restriction and interpolation. Cyclic reduction can be viewed as a two-level method with this Petrov-Galerkin coarse-grid operator and a block smoother (called F-relaxation) that converges in one iteration. As it is typical in the MGR setting, the MGRIT approach replaces the true Schur complement with a simpler operator (typically of the same form as the original bidiagonal system, but with a time propagator using time-step $m\Delta t$), replaces ideal restriction with simple injection, and compensates by adding relaxation. Furthermore, the two-level method can be extended to multiple levels in a simple recursive manner, and the full approximation storage (FAS) approach \cite{ABrandt_1977b} can be used to accommodate nonlinear problems. 

\section{Numerical results}
\label{sec:num_results}

We apply classical sequential time stepping and two MGRIT variants to the eddy current model problem \eqref{eq:pde} with the pulsed excitation \eqref{eq:input_pwm} on the space-time domain $\Omega \times (0,0.2]~\mathrm{s}$, with $\Omega = \Omega_0\cup\Omega_1\cup\Omega_2$ depicted in  Fig.~\ref{fig:model}. The spatial domain, $\Omega$, is discretized using $2269$ degrees of freedom and  the backward Euler method is used on a uniform grid with $32{,}768$ time steps for the time derivative of the space-discrete time-stepping problem \eqref{eq:magnetoquasi_discret}. The time step on the finest grid, $l=0$, is chosen to be $\Delta t = 6.1\cdot 10^{-6}\mathrm{s}$, and the time step on each coarse grid, $l$, is given by $m^l\Delta t$, $l\geq 1$. Two MGRIT variants are considered: a two-level Parareal-type method with a coarsening factor of $m=256$, and a five-level method that coarsens uniformly across all grids with a factor of $m=4$. Thus, the coarsest grid consists of $128$ points in time in both cases. On this coarsest temporal grid, time stepping is used. All spatial problems are solved using a direct LU solver.

The MGRIT algorithm was implemented in parallel using Python and Message Passing Interface (MPI). Numerical results were generated on an Intel Xeon Phi cluster consisting of 272 1.4 GHz Intel Xeon Phi processors.

\begin{figure}[!h]
	\begin{minipage}{.45\textwidth} 
		\includegraphics[width=\linewidth,height=40mm]{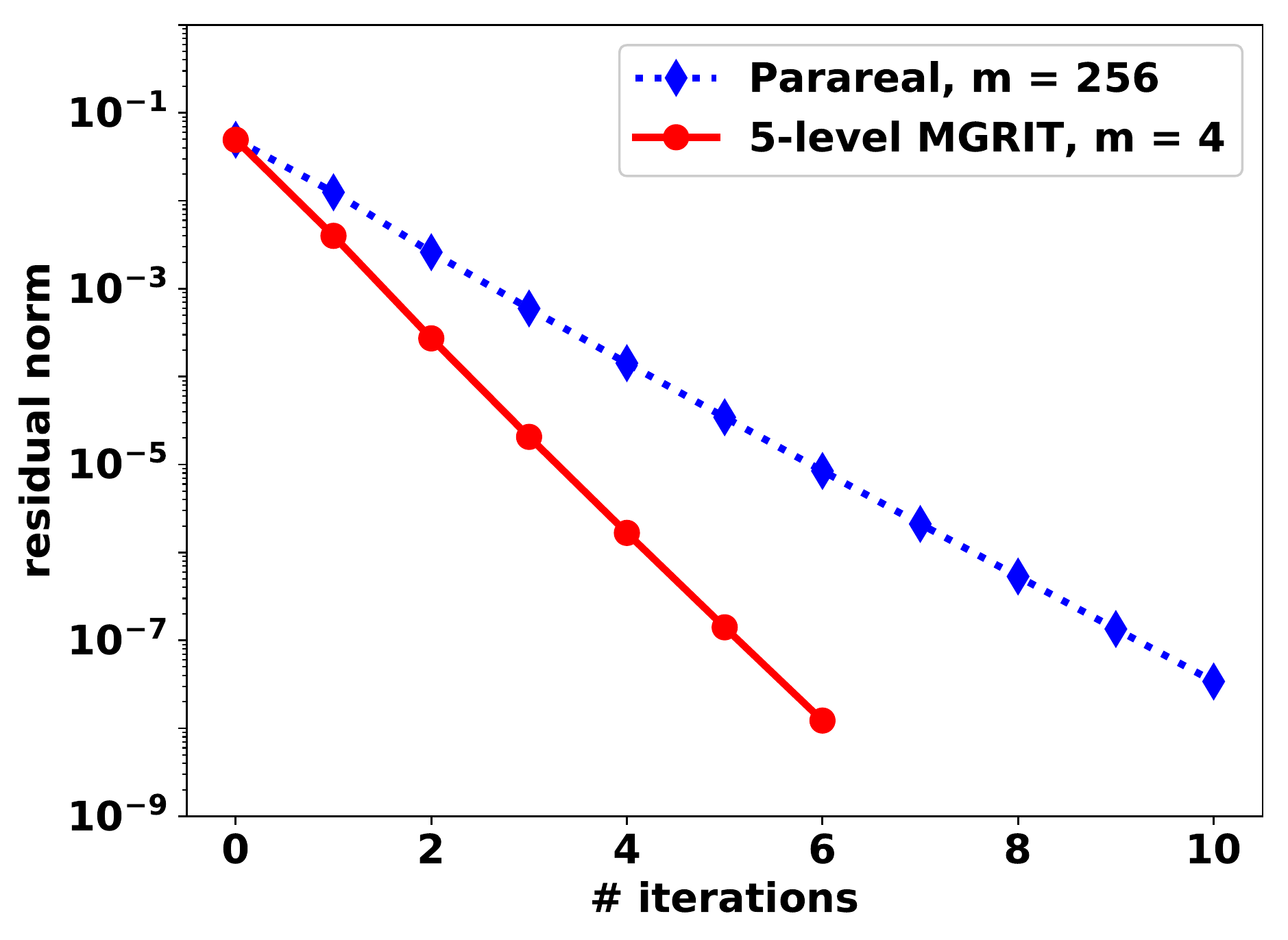} 
		\caption{Convergence of MGRIT variants applied to the eddy current model problem.} 
		\label{fig:convergence} 
	\end{minipage} \hfil \begin{minipage}{.45\textwidth} 
		\includegraphics[width=\linewidth,height=40mm]{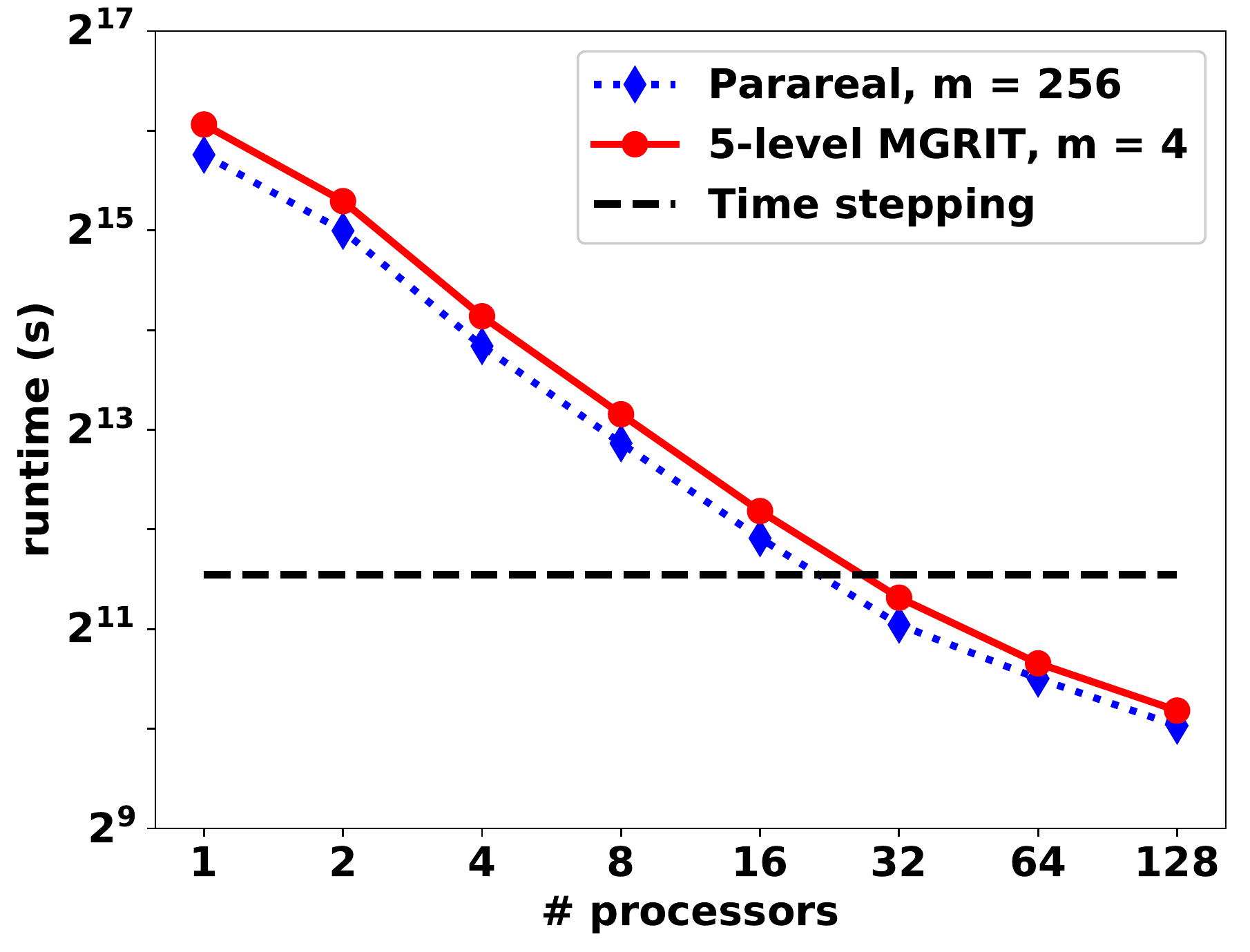} 
		\caption{Strong scaling results for MGRIT applied to the eddy current model problem.} 
		\label{fig:strong_scaling} 
	\end{minipage}
\end{figure}

Fig.~\ref{fig:convergence} shows convergence of the two MGRIT variants applied to the eddy current model problem. We see linear convergence for both variants. Comparing the number of spatial time-stepping solves required for the two methods to the optimal count of $N_t$ for sequential time stepping, we note that one iteration requires about $N_t$ or $2N_t$ spatial solves, respectively, when considering the two-level Parareal-type method or the five-level MGRIT scheme. This large computational overhead is demonstrated in the strong scaling results in Fig.~\ref{fig:strong_scaling}. The dotted and solid lines show results for the two- and five-level methods, respectively, for increasing the number of processors in the temporal dimension only. The dashed line shows the runtime of time stepping on one processor for reference purposes. Results show that the extra work in the MGRIT variants can be effectively parallelized at high processor counts, i.\,e., more than $32$, with good strong parallel scaling with a speedup of up to a factor of about $2.9$ over sequential time stepping.

\section{Conclusions}
\label{sec:conclusions}
MGRIT was applied for the first time to the eddy current problem, which yields an index-1 DAE after spatial discretization. A speedup of approximately three times could be obtained. A strong scaling investigation shows that the method converges linearly with the number of processors, even for non-standard, pulsed right-hand sides, which has been shown to be problematic for classical Parareal \cite{Kulchytska-Ruchka_2018ac}.

\begin{acknowledgement}
 The work is supported by the Excellence Initiative of the German Federal and State Governments, the Graduate School of Computational Engineering at TU Darmstadt, and the BMBF in the framework of project PASIROM (grants 05M18RDA and 05M18PXB).
\end{acknowledgement}

\end{document}